\documentclass[a4paper,12pt,dvipdfm]{amsart}
\usepackage[all,web]{xy}
\usepackage[dvips]{graphicx}
\usepackage{amsmath} 
\usepackage{amsfonts}
\usepackage{amssymb}
\usepackage{amsthm}
\usepackage{geometry}
\DeclareGraphicsExtensions{.eps}

\theoremstyle{plain}    
\newtheorem{thm}{Theorem}[section]
\numberwithin{equation}{section} 
\numberwithin{figure}{section} 
\theoremstyle{remark}   
 
\theoremstyle{plain}    
\theoremstyle{plain}    
 
\theoremstyle{remark} 

\theoremstyle{remark}

\theoremstyle{definition}

\theoremstyle{plain}
\newtheorem{conjecture}[thm]{Conjecture} 
\theoremstyle{plain}    
\theoremstyle{plain}    
\theoremstyle{definition}

\theoremstyle{definition}
\newtheorem{example}[thm]{Example}
\theoremstyle{plain}    
 
\theoremstyle{plain}    
\theoremstyle{remark}    
 
 \theoremstyle{remark}    
  
 \theoremstyle{definition}
  
 \theoremstyle{plain}
 \theoremstyle{remark}
 
 \theoremstyle{remark}    

\def\AA{\mathbb{A}}
\def\CC{\mathbb{C}}

\def\FF{\mathbb{F}}

\def\PP{\mathbb{P}}
\def\QQ{\mathbb{Q}}
\def\RR{\mathbb{R}}

\def\ZZ{\mathbb{Z}}

\font\cyr=wncyr10 scaled \magstep 1
\newcommand{\Sha}{\mbox{\cyr X}}

\newcommand{\mtwo}[4]{\left(
        \begin{matrix}#1&#2\\#3&#4
        \end{matrix}\right)}
\newcommand{\abcd}[4]{\left(
        \begin{smallmatrix}#1&#2\\#3&#4\end{smallmatrix}\right)}

\def\inner<#1,#2>{{\left\langle{{#1},{#2}}\right\rangle}}
\def\sl2of#1{\textrm{SL}_2(#1)}

\def\inner<#1,#2>{{\left\langle{{#1},{#2}}\right\rangle}}

\def\ssp{\def\baselinestretch{1.0}\large\normalsize}

\begin{document}

\title{On the index of the Heegner subgroup of elliptic curves}
\author{Carlos Casta\~no-Bernard}
\urladdr{http://users.ictp.it/~ccastano/}
\email{ccastano@ictp.it}

\begin{abstract}
Let $E$ be an elliptic curve of conductor $N$ and rank one over $\QQ$.
So there is a non-constant morphism $X_0^+(N)\longrightarrow E$
defined over $\QQ$,
where $X_0^+(N)=X_0(N)/w_N$ and $w_N$ is the Fricke involution.
Under this morphism
the traces of the Heegner points of $X_0^+(N)$
map to rational points on $E$.
In this paper we study
the index $I$ of the subgroup generated by all these traces on $E(\QQ)$.
We propose and also discuss a conjecture that says that
if $N$ is prime and $I>1$,
then either the number of connected components $\nu_N$ of 
the real locus $X_0^+(N)(\RR)$ is $\nu_N>1$
or (less likely) the order $S$ of the Tate-\v Safarevi\v c group $\Sha(E)$ of $E$ is $S>1$.
This conjecture is backed by computations performed on each $E$
that satisfies the above hypothesis in the range $N\leq 129999$.

This paper was prepared for the proceedings of
the Conference on
Algorithmic Number Theory,
Turku, May 8--11, 2007.
We tried to make the paper as self contained as possible.
\end{abstract}

\address{Mathematics Section, ICTP, Strada Costiera 11, I-34014 Trieste (Italy)}

\maketitle

\pagenumbering{roman}
\setcounter{page}{0}



\tableofcontents


\pagenumbering{arabic}
\pagestyle{headings}

%
%
%
\section{Introduction}

%
\subsection{Motivation}
Let $E$ be an elliptic curve over $\QQ$,
i.e. a complete curve of genus one
with a specified rational point $O_E$,
hence $E$ has a natural structure of a commutative algebraic group
with zero element $O_E$.
The Mordell-Weil theorem asserts that
the group $E(\QQ)$ of rational points on $E$
is finitely generated.
So the classical Diophantine problem
of determining $E(\QQ)$
is thus the problem of obtaining a finite set of generators for
the group $E(\QQ)$.
The finite subgroup $E(\QQ)^\textit{tors}$
of torsion points of $E(\QQ)$ is easy to compute.
However,
finding generators $g_1,\dots,g_{r_E}$
for the free abelian group $E(\QQ)/E(\QQ)^\textit{tors}$
is in general a hard problem.
The Birch and Swinnerton-Dyer conjecture predicts
(among other things)
that the rank $r_E$ of
the Mordell-Weil group $E(\QQ)$ is
the order of vanishing at $s=1$ of the Hasse-Weil $L$-function $L(E,s)$
attached to $E$.
By the work of Kolyvagin on Euler systems of Heegner points on
(certain twists of) modular elliptic curves,
and the well-known fact due to
Wiles~\cite{wiles:flt},
and Breuil-Conrad-Diamond-Taylor~\cite{breuil:modularity}
that every elliptic curve $E$ over $\QQ$ admits
a (non-constant) morphism $\varphi:X_0(N)\longrightarrow E$ over $\QQ$,
we know that this prediction is true for $r_E=0$ and $1$.
We are interested in the latter case,
and henceforth we assume that $L(E,s)$ has a simple zero at $s=1$.
Then $\varphi$ factors through
the quotient $X_0^+(N)=X_0(N)/w_N$ associated to the Fricke involution $w_N$
and the so-called Heegner point construction\footnote{Heegner points
were first studied systematically by Birch~\cite{birch:heegner}.}
yields a non-trivial subgroup $H$ of $E(\QQ)/E(\QQ)^\textit{tors}$.
Gross-Kohnen-Zagier~\cite[p. 561]{gross:gkz} proved that
the (full) Birch and Swinnerton-Dyer for $r_E=1$ is equivalent to
\begin{equation}\label{eqn:gkz}
I_E^2 = c_E\cdot n_E\cdot m_E\cdot |\Sha(E)|,
\end{equation}
where $I_E$ is the index of $H$, $\Sha(E)$ is
the Tate-\v Safarevi\v c group of $E$, $c_E$ is Manin's constant, $m_E$ is the product of
the Tamagawa numbers,
and $n_E$ is the index of a certain subgroup of
the $-1$-eigenspace $H_1(E(\CC);\ZZ)^{-}$ of
complex conjugation acting on $H_1(E(\CC);\ZZ)$
constructed in terms of classes of
Heegner geodesic cycles in $H_1(X_0^+(N)(\CC);\ZZ)^{-}$.
(The relevant definitions are recalled below.)
Let us assume this conjecture.
To simplify our discussion let us assume further that
the conductor $N_E$ of $E$ is prime so that
the index $I_E$ is completely determined by $n_E$  and $|\Sha(E)|$.
Numerical evidence strongly suggests that
there are $109$ curves such that $I_E>1$ out of
the  $914$ curves $E$ of rank one and prime conductor $N\leq 129999$
in Cremona's Tables~\cite{cremona:onlinetables}.
For each of these curves with $I_E>1$,
then either the number $\nu_N$ of connected components of 
the real locus $X_0^+(N)(\RR)$ of the quotient modular curve $X_0^+(N)$
is $\nu_N > 1$ or, less likely (only $8$ cases), $\Sha(E)$ is non-trivial.
This suggests a non-trivial connection between the topology of $X_0^+(N)(\RR)$
and the arithmetic of $E$,
which is not expected since $\nu_N$ is a certain simple sum of
class numbers of real quadratic fields and
heuristic considerations suggest that
the equality $\nu_N = 1$ is more likely than
the inequality $\nu_N > 1$.
This paper is about a conjecture 
motivated by the above discussion.
We state it in Subsection~\ref{sub:dis}
and then discuss a homological formulation of our conjecture
which hopefully will furnish a new approach to Equation~\ref{eqn:gkz}.

%
\subsection{Acknowledgements}
I would like to heartily thank Professor Birch,
whose comments encouraged me to investigate further
some ``loose ends'' related to some odd behaviour for the curve \textbf{359A}
mentioned in my Ph.\thinspace{}D. thesis~\cite[p. 75]{ccb:thesis}.
I would also like to thank my colleagues
at ICTP whose financial support,
through the granting of a visiting fellowship,
has facilitated the writing of this paper.

Table~\ref{tbl:indexes01} and Table~\ref{tbl:indexes02}
were computed with the help of \textsc{Pari}~\cite{pari:gp},
installed on GNU/Linux computers.

%
%
%
\section{Background}\label{sec:bac}

\subsection{The Hasse principle and genus one curves}
It is a classical Diophantine problem
the determination of the set of rational points $C(\QQ)$ of a given
complete non-singular algebraic curve defined over $\QQ$.
The problem is solved for the case of genus zero.
Legendre theorem,
as stated by Hasse,
says that given any conic $C$ with coefficients in $\QQ$
the set $C(\QQ)$ is non-empty if and only if
the set $C(\QQ_p)$ is non-empty
for every prime $p$ including $p=\infty$,
where $\QQ_p$ is the field of $p$-adic numbers,
if $p\not=\infty$ and $\QQ_p=\RR$,
if $p=\infty$. 
Moreover,
it is known that it suffices to determine whether $C(\QQ_p)$ is non-empty
for each prime $p$ that divides the discriminant $D$ of
an homogeneous equation $f(X,Y,Z)=0$ for the conic $C$.
Then by Hensel's lemma we know that $f(X,Y,Z)=0$
will have a non-trivial zero in $\QQ_p$ for $p|D$ if and only if it has
an ``approximate'' zero.
Once we have a rational point $O$ on $C$,
it is easy to see that there are an infinite number of them
by fixing any line $L\subset\PP^2$ defined over $\QQ$ (e.g. the $X$-axis)
and parametrise $C(\QQ)$ with $L$ in the obvious way.
This furnishes an algorithm to effectively compute $C(\QQ)$
in the genus zero case.

Let us consider the genus one case.
By the work of Selmer~\cite{selmer:dioph01} we know that
the obvious extension of Legendre's theorem to curves of genus one is not true.
For example the curve $C$ in $\PP^2$ given by the Selmer cubic
\begin{displaymath}
3X^3 + 4Y^3 +5Z^3=0
\end{displaymath}
is such that $C(\QQ_p)\not=0$ for every prime $p$,
including $p=\infty$.
But it turns out that $C(\QQ)=\emptyset$.
In such cases it is said that $C$ violates the \textit{Hasse principle}.
There is a natural way to measure the extent of failure of this principle.
The \textit{Jacobian} $E=\textit{Jac}(C)$ of $C$
is a complete non-singular genus one curve defined over $\QQ$ equipped with
a commutative algebraic group structure,
i.e. $E$ is an elliptic curve,
together with an isomorphism $j:C\longrightarrow E$ over $\QQ^\textit{alg}$
such that for every element $\sigma$ in the Galois group $G_\QQ$
of $\QQ^\textit{alg}$ over $\QQ$ the map
\begin{displaymath}
(\sigma\circ j)\circ j^{-1}:\textit{Jac}\,(T)\longrightarrow\textit{Jac}\,(T)
\end{displaymath}
is of the form  $P\mapsto P + a_\sigma$, for some $a_\sigma\in E(\QQ^\textit{alg})$.
So we may define
the \textit{Tate-\v Safarevi\v c group} $\Sha(E)$ of $E$ as
the set of isomorphism classes of pairs $(T,\iota)$,
where $T$ is a smooth curve defined over $\QQ$ of genus one
such that $T(\QQ_p)\not=\emptyset$,
for all $p$ prime and $\iota:E\longrightarrow\textit{Jac}\,(T)$
is an isomorphism defined over $\QQ$.
(Given $T$ such that $E=\textit{Jac}(T)$,
the map $\sigma\mapsto a_\sigma$ is a $1$-cocycle
whose image in the cohomology group $H^1(G_\QQ, E)$ is
uniquely determined by the isomorphism class of $(T,\iota)$.
So we may identify $\Sha(E)$ with a subgroup of $H^1(G_\QQ, E)$.
Cf. Cassels book~\cite{cassels:elliptic}.)
Clearly the Hasse principle holds for $C$
if and only if $\Sha(E)$ consists of exactly one element,
where $E$ is the Jacobian of $C$.
It is conjectured to be finite,
i.e. that Hasse principle fails by a ``finite amount'' in the genus one case.

Cassels' proved that if $\Sha(E)$ is indeed finite,
then its order is a square.

%
\subsection{Structure of the Mordell-Weil group}
\begin{figure}
\begin{center}
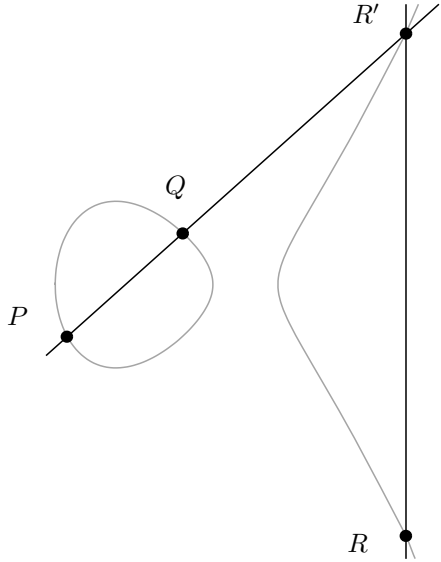 
\end{center}
\caption{Group law $P+Q=R$.}
\label{cap:chordtang}
\end{figure}

The algebraic group structure of an elliptic curve $E$ may be made explicit as follows.
Let $O_E$ be the zero element of $E$.
Using the Riemann-Roch theorem we see that
the map Albanese map $P\longmapsto P-O_E$
identifies the set $E(K)$ of $K$-rational points of $E$ with
the Picard group $\textrm{Pic}^0(E/K)$ of $E$ over any field $K$ containing $\QQ$.
Using again the Riemann-Roch theorem we may see that $E$ has
a \textit{Weierstra\ss\ model}
\begin{equation}\label{eqn:w}
Y^2Z+a_1XYZ+a_3YZ^2=X^3+a_2X^2Z+a_4XZ^2+a_6Z^3
\end{equation}
where $O_E$ corresponds to $(0:1:0)$,
for $a_1$, $a_2$, $a_3$,$a_4$, and $a_6\in \QQ$
such that the discriminant $\Delta$ of Equation~\ref{eqn:w} is non-zero.
It is well-known that the converse holds,
so a curve defined by a Weierstra\ss\ equation such that $\Delta\not=0$
is a complete non-singular curve of genus one,
and thus an elliptic curve with zero element $O_E=(0:1:0)$.
In particular,
the curve obtained by reducing the coefficients of Equation~\ref{eqn:w}
modulo a prime number $p$
is an elliptic curve if and only if $p$ does not divide $\Delta$,
in which case we say that $E$ has \textit{good reduction} at $p$.
A further consequence of the Riemann-Roch theorem is that
the group law is given by the classical chord and tangent construction,
which is schematically outlined in Figure~\ref{cap:chordtang}.
Using this geometric property we may easily write down
explicit rational functions with coefficients in $\QQ$ on
the coordinate functions $x$ and $y$ for 
the addition law $E\times E\longrightarrow E$
and for the inverse of an element law $E\longrightarrow E$.

The Mordell-Weil theorem asserts that the group $E(\QQ)$
is a finitely generated abelian group,
thus $E(\QQ)=E(\QQ)^\textit{tors}\oplus E(\QQ)^\textit{free}$,
where
the torsion subgroup $E(\QQ)^\textit{tors}\subset E(\QQ)$ is finite
and $E(\QQ)^\textit{free}\subset E(\QQ)$ is a free subgroup of (finite) rank $r_E$.
It is well-known that
the subgroup $E(\QQ)^\textit{tors}$ is not difficult to compute.
However,
obtaining generators for a subgroup $E(\QQ)^\textit{free}$ is in general a hard problem.
A measure of the arithmetic complexity of a given non-torsion rational point $P$ on $E$
is given by its \textit{N\'eron-Tate height}
\begin{displaymath}
\hat{h}(P)=\lim_{n\rightarrow\infty}4^{-n}h(2^nP),
\end{displaymath}
where the \textit{na\"\i ve height} $h(P)$ of a point $P=(x:y:z)$ in $\PP^2(\QQ)$
is given by $h(P)=\log\max(|x|,|y|,|z|)$,
where $x$, $y$, and $z$ are integers such that $\gcd(x,y,z)=1$.
It is well-known that $\hat{h}(P)$
does not depend on the choice of Weierstra\ss\ model for $E$ and,
moreover,
it defines
a non-degenerate positive definite quadratic form
on the $r_E$-dimensional real vector space $E(\QQ)\otimes_\ZZ\RR$.
The \textit{height paring} is
the bilinear form $\langle\cdot,\cdot\rangle$ on $E(\QQ)\otimes_\ZZ\RR$
such that $\langle P,P\rangle=\hat{h}(P)$,
for all $P\in E(\QQ)\otimes_\ZZ\RR$.
The determinant $R_E$ of the $r_E$ by $r_E$ matrix whose entries are given by
the height paring $\langle\cdot,\cdot\rangle$ applied to
a set of generators of $E(\QQ)^\textit{free}$
is known as the \textit{regulator} of $E(\QQ)$.

%
\subsection{The Birch and Swinnerton-Dyer conjecture}
As above let $E$ be an elliptic curve defined over $\QQ$,
and suppose we have used Tate's algorithm~\cite{tate:algorithm}
to obtain the conductor $N_E$
and a \textit{minimal Weierstra\ss\ model} of $E$,
i.e. an integral Weierstra\ss\ model of $E$
with $|\Delta|$ minimal.
Such discriminant is known as
the \textit{minimal discriminant}\footnote{The minimal discriminant $\Delta_E$
and the conductor $N_E$ share the same prime divisors,
and under certain circumstances they coincide (up to multiplication by $\pm 1$),
e.g. when $\Delta_E$ is prime.}
of $E$ and denote we it $\Delta_E$.
The Hasse-Weil $L$-function of $E$ over $\QQ$ is
\begin{displaymath}
L(E,s)=\sum_{n=1}^\infty a_E(n)n^{-s}=
\prod_{\textit{prime}\,p\mid N}
(1-a_E(p)p^{-s})^{-1}
\prod_{\textit{prime}\,p\nmid N}
(1-a_E(p)p^{-s} + p^{1-2s})^{-1},
\end{displaymath}
where
\begin{displaymath}
a_E(p)=\left\{
\begin{array}{ll}
p+1-\#(E(\FF_p)),&\textit{good reduction,}\\
1,&\textit{split reduction,}\\
-1,&\textit{non-split reduction,}\\
0,&\textit{cuspidal reduction.}\\
\end{array}
\right.
\end{displaymath}
Since $E$ is defined over $\QQ$ the work of
Wiles~\cite{wiles:flt}
and
Breuil-Conrad-Diamond-Taylor~\cite{breuil:modularity}
implies that $E$ is modular,
and in particular $L(E,s)$ may be analytically continued to the whole 
complex plane $\CC$.
(See below.)
The Birch and Swinnerton-Dyer conjecture predicts that $L(E,s)$
has a Taylor expansion around $s=1$ of the form
\begin{displaymath}
L(E,s)=\kappa_{r_E}(s-1)^{r_E}+\kappa_{r_E+1}(s-1)^{r_E+1}+\dots,
\end{displaymath}
where
\begin{displaymath}
\kappa_{r_E}= |\Sha(E)|m_E
\frac{R_E}{|E(\QQ)^\textrm{tors}|}\, \Omega_E
\end{displaymath}
where $m_E$ is the product of all the local Tamagawa numbers $c_p$,
and $\Omega_E$ is the least positive real period of the N\'eron differential
\begin{equation*}
\omega_E=\frac{dX}{2Y+a_1X+a_3},
\end{equation*}
where $a_1$ and $a_3$ are as in Equation~\ref{eqn:w}
(assuming the Weierstra\ss\ equation is minimal).

\begin{example}
The Selmer cubic $C$ defined by $3X^3 + 4Y^3 + 5Z^3 = 0$
has Jacobian $E$ with Weierstra\ss\ model $Y^2 = 4X^3 - 97200$.
(Cf. Perlis~\cite[p. 58]{perlis:thesis}.)
Using Tate's algorithm we may see that $E$ has conductor $N_E=24300$
and minimal Weierstra\ss\ model $Y^2 = 4X^3 - 24300$.
Using this information we may identify $E$ 
in entry \textbf{24300 Y 2} of Cremona's Tables~\cite{cremona:onlinetables}.
According to that entry the rank of $E$ is zero and
the order of Tate-\v Safarevi\v c group predicted by
the Birch and Swinnerton-Dyer conjecture is $|\Sha(E)|=3^2$.
This is consistent with the fact that the Hasse principle fails for $C$,
as remarked above.
\end{example}


%
%
%
\section{On the index $I_\varphi$ and the topology of $X_0^+(N)(\RR)$}\label{sec:ind}

%
\subsection{Modular parametrisation}
Let $X_0(N)$ be the normalisation of
the moduli space that classifies
pairs $(A,A^\prime)$ of elliptic curves
together an isogeny $\phi:A\longrightarrow A^\prime$
with cyclic kernel of order $N$.
The curve $X_0(N)$ may be described as follows.
Let $\Gamma$ be
the group $\textrm{SL}_2(\RR)=\{\abcd{a}{b}{c}{d}:ad-cb=1\}$
modulo multiplication by $\pm 1$,
and let $\Gamma$ act on
the upper half plane $\mathfrak{h}=\{z\in\CC:\Im(\tau)>0\}$ in the usual way by letting
\begin{displaymath}
\tau\mapsto\frac{a\tau+b}{c\tau+d}.
\end{displaymath}
First we may identify the complex points of
the moduli space $Y(1)$ that classifies elliptic curves $E$ over $\CC$ with 
the complex points of the affine line $\AA^1$
by mapping the isomorphism class of $E\cong \CC/(\ZZ\tau + \ZZ)$
to the image of $\tau$ in  $\Gamma\backslash\mathfrak{h}$
followed by the classical $j$-invariant map
\begin{displaymath}
j(\tau) = \frac{E_4^3}{\Delta}(\tau) =\frac{1}{q} + 744 + 196884q + \dots,
\end{displaymath}
where $\Delta$ is the cusp form
of weight $12$ defined by the infinite product $\Delta(\tau)=q\prod_{n>0}(1-q^n)^{24}$,
and $E_4$ is the modular form of weight $4$ defined by
the series $E_4(\tau)=1+240\sum_{n>0}\sigma_3(n)q^n$,
where as usual $\sigma_k(n)=\sum_{0<d|n}d^k$ and $q=e^{2\pi i\tau}$.
The obvious action of $\Gamma$ on
the \textit{cusps} $\PP^1(\QQ) = \QQ\cup\{i\infty\}$
is transitive,
so the (one-point) compactification $X(1)(\CC)$ of
the complex line $Y(1)(\CC)$ is the Riemann sphere $X(1)=\Gamma\backslash\mathfrak{h}^*$,
where $\mathfrak{h}^*=\mathfrak{h}\cup\PP^1(\QQ)$.
We also have a bijection
\begin{displaymath}
\xymatrix{
\Gamma_0(N)\backslash\mathfrak{h}^*\ar [r]& X_0(N)(\CC)\\
\tau\mod{\Gamma_0(N)}\ar @{|->}[r]& [\CC/(\ZZ\tau + \ZZ)\longrightarrow\CC/(\ZZ\tau + \frac{1}{N}\ZZ)]\\
}
\end{displaymath}
where 
\begin{equation*}
\Gamma_0(N)=
\left\{\mu=\mtwo{a}{b}{c}{d}\in\textrm{SL}_2(\ZZ)
\,\colon\,
c\equiv 0\pmod{N}\right\}.
\end{equation*}
The quotient set $\Gamma_0(N)\backslash\mathfrak{h}^*$
has a unique complex-analytic structure such that the
natural map $\psi:\Gamma_0(N)\backslash\mathfrak{h}^*\longrightarrow X(1)(\CC)$
is a proper.
Moreover,
the above bijection is in fact an isomorphism 
between $\Gamma_0(N)\backslash\mathfrak{h}^*$ and $X_0(N)(\CC)$
as Riemann surfaces in such a way that $\psi$ is induced by
the projection map $(A,A^\prime)\mapsto A$.
The degree of $\psi$ is the degree of
the minimum polynomial $\Phi_N(j,Y)\in \CC(j)[Y]$
of $j(N\tau)$ over $\CC(j)$,
and it turns out that $\Phi_N(X,Y)$ has integral coefficients.
The field of fractions of $\QQ[X,Y]/(\Phi_N(X,Y))$ gives the
canonical $\QQ$-structure of $X_0(N)$.

The \textit{Fricke involution} $w_N$
may be defined as the morphism of $X_0(N)$ to itself
induced by
mapping an isogeny $\phi:A\longrightarrow A^\prime$
to its dual $\hat{\phi}:A^\prime\longrightarrow A$.
In the complex-analytic setting $w_N$ is induced by
the involution $\tau\mapsto-\frac{1}{N\tau}$ of $\mathfrak{h}$.
Let $X_0^+(N)$ be the quotient of $X_0(N)$ by the group $\{1,w_N\}$.
The classical result $\Phi_N(X,Y)=\Phi_N(Y,X)$ implies that
the canonical map  $X_0(N)\longrightarrow X_0^+(N)$ is defined over $\QQ$.

Again let $E$ be an elliptic curve defined over $\QQ$.
As mentioned above,
by the work of
Wiles~\cite{wiles:flt}
and
Breuil-Conrad-Diamond-Taylor~\cite{breuil:modularity}
we know that $E$ is modular.
This means that the Fourier series $f_E(\tau)=\sum_n a_E(n)q^n$
is a normalised newform,
and thus $\omega_f=2\pi i f_E(\tau)d\tau$
is a holomorphic differential
on $X_0(N)$ such that the map $\varphi:X_0(N)\longrightarrow E$
defined by
\begin{equation*}
\tau\mod{\Gamma_0(N)}\mapsto\int_{i\infty}^\tau \omega_f
\end{equation*}
followed by
the classical map $z\mapsto (\wp_\Lambda(z),\wp_\Lambda^\prime(z))$,
is a (well-defined) non-constant morphism over $\QQ$,
where $\wp_\Lambda$ is the Weierstra\ss\ $\wp$-function
and $\Lambda\subset\CC$ is the lattice generated by
the periods of a N\'eron differential $\omega_E$
associated to a minimal Weierstra\ss\ model of $E$.
From now on we assume that $E$ has rank one over $\QQ$.
By the work of Kolyvagin~\cite{kolyvagin:euler},
Gross-Kohnen-Zagier~\cite{gross:gkz}
and results due to 
Waldspurger, Bump, Friedberg and Hoffstein,
we know that 
if $r=0$ or $1$,
then the order of vanishing of $L(E,s)$ is as predicted by
the Birch and Swinnerton-Dyer conjecture
(and also that $\Sha(E)$ is finite).
In particular $L(E,s)$ has a simple zero at $s=1$ and thus $w_N\omega_f=\omega_f$.
So the modular parametrisation factors through
the quotient map $X_0^+(N)\longrightarrow X_0(N)$.

%
\subsection{Heegner points}
Now suppose we fix a pair of integers $(D,r)$
that satisfy the so-called 
\textit{Heegner condition}\footnote{This condition was 
introduced by Birch~\cite{birch:heegner}.}
i.e. $D$ is the discriminant of
an imaginary quadratic order $\mathcal{O}_D$
of conductor $f$ such that gcd$(N,f)=1$
and $r\in\ZZ$ is such that
\begin{equation*}
D\equiv r^2\pmod{4N}.
\end{equation*}
So we have
a proper $\mathcal{O}_D$-ideal $\mathfrak{n}_r=\ZZ N + \ZZ\frac{-r+\sqrt{D}}{2}\subset K=\QQ(\sqrt{D})$
and $\mathcal{O}_D/\mathfrak{n}_r\cong \ZZ/N\ZZ$.
So for each proper $\mathcal{O}_D$-ideal $\mathfrak{a}\subset K$ we have
a point $x=(\CC/\mathfrak{a},\CC/(\mathfrak{n}_r^{-1}\mathfrak{a}))$
on $X_0(N)$.
This point $x$ is known as a \textit{Heegner point},
and following Gross~\cite{gross:heegner}
we denote it by $x=(\mathcal{O}_D,\mathfrak{n}_r,[\mathfrak{a}])$,
where $[\mathfrak{a}]$ is the 
class of $\mathfrak{a}$ in Pic$(\mathcal{O}_D)$.
The latter set may be identified with the $\Gamma$-orbits $\Gamma\backslash\mathcal{Q}_D^0$
of the set $\mathcal{Q}_D^0$ of
primitive binary quadratic forms $[A,B,C]$ of discriminant $D=B^2-4AC$ and $A>0$
by writing each $\mathcal{O}_D$-ideal $\mathfrak{a}$
as $\mathfrak{a}=A\ZZ + \frac{-B+\sqrt{D}}{2}\ZZ$,
for some $[A,B,C]\in\mathcal{Q}_D^0$.
Moreover,
the $\Gamma_O(N)$-orbits $\Gamma_O(N)\backslash\mathcal{Q}_{N,D,r}^0$
of the set $\mathcal{Q}_{N,D,r}^0$ of $[A,B,C]\in\mathcal{Q}_D^0$
such that $N|A$ and $B\equiv r\pmod{2N}$ may be identified with
the set of Heegner points $(\CC/\mathfrak{a},\CC/(\mathfrak{n}_r^{-1}\mathfrak{a}))$,
and also with
the set of $\Gamma_O(N)$-orbits of points $\tau\in\mathfrak{h}$
of the form $\tau=\frac{-B+\sqrt{D}}{2A}$.

The field of definition $H$ of each Heegner point $x=(A,A^\prime)$
may be described as follows.
Note that a point $x=(A,A^\prime)$ on $X_0(N)(\CC)$ is a Heegner point
associated to $D$
if and only if $\textit{End}(A)=\textit{End}(A^\prime)=\mathcal{O}_D$.
So $H=K(j(\tau))$ where $\tau=\frac{-B+\sqrt{D}}{2A}$ is as above,
and by the theory of Complex Multiplication
the action of the Galois group $\textrm{Gal}(K^\textit{alg}/K)$
on $x$ is determined by a homomorphism 
\begin{equation*}
\delta\colon\textrm{Gal}(K^\textit{alg}/K)\longrightarrow\textrm{Pic}(\mathcal{O}_D)
\end{equation*}
such that $\delta(\sigma)*x = x^\sigma$,
where $*$
is defined by $\mathfrak{b}*x=(\mathcal{O}_D,\mathfrak{n}_r,[\mathfrak{b}^{-1}\mathfrak{a}])$.
In other words $H$ is
the fixed field of the Galois group $\textrm{ker}(\delta)$
and $\textrm{Gal}(H/K)\cong\textrm{Pic}(\mathcal{O}_D)$.
The field $H$ is known as the \textit{ring class field} attached to $\mathcal{O}_D$,
i.e. the maximal abelian extension of $K$
unramified at all primes $\mathfrak{p}$ of $K$ which do not divide $f$.
More precisely,
the homomorphism $\delta$ is the inverse of the Artin reciprocity map,
so in fact $\delta(\textrm{Frob}_{\mathfrak{p}})=[\mathfrak{p}]$
for each prime $\mathfrak{p}$ of $K$ which does not divide $f$,
where $\textrm{Frob}_{\mathfrak{p}}\in\textrm{Gal}(H/K)$ is
the Frobenius element at $\mathfrak{p}$,
which is characterised by the properties $\textrm{Frob}_{\mathfrak{p}}\mathfrak{P}=\mathfrak{P}$
and $\textrm{Frob}_{\mathfrak{p}}\alpha\equiv\alpha^q\pmod{\mathfrak{P}}$,
for each $\alpha$ in the ring of integers $\mathcal{O}_H$ of $H$,
where $\mathfrak{P}$ is a prime ideal of $H$ above $\mathfrak{p}$
and $q=\#(\mathcal{O}_K/\mathfrak{p})$.

To simplify the exposition we assume from now on that
the discriminant $D$ is fundamental,
and also that $E(\QQ)\cong \ZZ$.
The \textit{weighted trace} $y_{D,r,\varphi}$ on $E$
associated to the pair $(D,r)$ may be defined by the equation
\begin{equation}\label{eqn:trace}
u_Dy_{D,r,\varphi}=
\sum_{\mathfrak{a}\in\textrm{Pic}(\mathcal{O}_D)} \varphi(\mathcal{O}_D,\mathfrak{n}_r,[\mathfrak{a}]),
\end{equation}
where
\begin{displaymath}
u_D=\left\{
\begin{array}{ll}
\frac{1}{2}\#(O^\times_D),&\textrm{if $\#(O^\times_D)>2$.}\\ 
\\
2,&\textrm{if $\#(O^\times_D)=2$ and $N|D$,}\\
\\
1,&\textrm{otherwise.}\\
\end{array}
\right.
\end{displaymath}
We claim that $y_{D,r,\varphi}$ is a rational point on $E$.
Since $K$ is an imaginary quadratic field,
the non-trivial element of $\textrm{Gal}(K/\QQ)$
is complex conjugation,
which acts on Heegner points
as $(\mathcal{O}_D,\mathfrak{n}_{r},[\mathfrak{a}])\mapsto(\mathcal{O}_D,\mathfrak{n}_{-r},[\mathfrak{a}^{-1}])$.
Also,
note that the action of the Fricke involution $w_N$ is given
by $w_N(\mathcal{O}_D,\mathfrak{n}_r,[\mathfrak{a}])=(\mathcal{O}_D,\mathfrak{n}_{-r},[\mathfrak{n}^{-1}\mathfrak{a}])$.
Therefore the action of $w_N$ on
the right-hand side of Equation~\ref{eqn:trace}
is the same as that of complex conjugation.
But we assumed $\varphi$ factors through
the canonical quotient map $X_0(N)\longrightarrow X_0^+(N)$
associated to $w_N$.
Thus the right-hand side of Equation~\ref{eqn:trace}
is defined over $\QQ$.
Finally,
each Heegner point $\tau\in\mathfrak{h}$ of discriminant $D$
is the fixed point of an element of order $u_D$ of the group generated by $\Gamma_0(N)$
and the Fricke involution $w_N$ (cf. Zagier~\cite{zagier:modular}),
and our claim follows.

Recall we assumed $E(\QQ)\cong \ZZ$.
So we may fix a generator $g_E$ of
the Mordell-Weil group $E(\QQ)$ of $E$ over $\QQ$.
The index $I_{D,r,\varphi}$ of $y_{D,r,\varphi}$ in $E(\QQ)$
may be expressed as
\begin{equation*}
y_{D,r,\varphi}=I_{D,r,\varphi}\,g_E,
\end{equation*}
We are interested in the index $I_\varphi$ of
the group generated by the Heegner points,
i.e. the greatest common divisor of the indexes $I_{D,r,\varphi}$
for all pairs $(D,r)$ that satisfy the Heegner condition.

%
\subsection{Heegner paths}
Suppose that the pair $(\Delta,\rho)$
satisfies the Heegner condition.
Suppose further that $\Delta>0$
and that $\Delta$ is not the square of an integer.
Assume the above notation
and let $Q=[A, B, C]\in\mathcal{Q}_{N,\Delta,\rho}^0$
The condition $N|A$ implies that
all the automorphs of $Q$ lie in $\Gamma_0(N)$.
More explicitly,
if $(x,y)\in\ZZ\times\ZZ$ is a fundamental solution of
Pell's equation $X^2 - DY^2 = 1$ then
the fundamental automorph of $Q$ given by
\begin{equation*}
M_Q = \mtwo{x - By}{-2Cy}{2Ay}{x + By}
\end{equation*}
lies in $\Gamma_0(N)$.
Note that $M_Q$ fixes $\tau^\pm(Q)=\frac{-B\pm\sqrt{\Delta}}{2A}\in\PP^1(\RR)$.
We normalise our choice of $M_Q$ by assuming that
the eigenvalue $\lambda_Q=x + y\sqrt{\Delta}\in\mathcal{O}_D^\times$
is $\lambda_Q>1$,
so that $\tau^-(Q)$ is repelling and $\tau^+(Q)$ is attracting.
The \textit{axis} of $M_Q$ is
the geodesic $\{\tau^-(Q),\tau^+(Q)\}\subset\mathfrak{h}$
from $\tau^-(Q)$ to $\tau^+(Q)$.
Clearly it is stable under the action of $M_Q$
and has the same orientation as
the geodesic segment $\{\tau_0,M_Q\tau_0\}$,
given any point $\tau_0$ on it.
Now let $\gamma_{Q,\tau_0}$ be
the closed path on $X_0(N)(\CC)$ defined by $\{\tau_0,M_Q\tau_0\}$.
It is a smooth path on $X_0(N)(\CC)$
except when it contains an elliptic point of order $2$,
in which case $\gamma_{Q,\tau_0}=-\gamma_{Q,\tau_0}$ as $1$-cycles.
Note $\gamma_{Q,\tau_0}$ depends only on
the $\Gamma_0(N)$-equivalence class of $Q$ 
So given $(D_0,r_0)$ and $(D_1,r_1)$ that satisfy
the Heegner condition
we may define the \textit{(twisted) Heegner cycle}
\begin{equation*}
\gamma_{D_0,D_1,\rho}=\sum_{[Q]\in\Gamma_0(N)\backslash\mathcal{Q}_{N,\Delta,\rho}^0}\chi_{D_0}(Q)\gamma_Q
\end{equation*}
where $\chi_{D_0}$ is the \textit{generalised genus character},
following Gross-Kohnen-Zagier~\cite[p. 508]{gross:gkz};
\begin{displaymath}
\chi_{D_0}(Q)=\left\{
\begin{array}{ll}
\left(\frac{D_0}{n}\right),&\textrm{if $\textrm{gcd}(A/N,B,C,D_0)=1$}\\ 
&\\ 
0,&\textrm{otherwise.}\\
\end{array}
\right.
\end{displaymath}
where $\Delta=D_0D_1$ and $\rho=r_0r_1$.
Here in the first case $n$ is an integer represented by $[A/N^\prime,B,CN^\prime]$,
where $N^\prime$ is a positive divisor of $N$,
and $Q=[A,B,C]$.
Note $\gamma_{D_0,D_1,\rho}$ is invariant with respect the action of
the Fricke involution $w_N$,
so it defines a $1$-cycle on the quotient Riemann surface $X_0^+(N)(\CC)$.
If we assume further that $D_0<0$ and $D_1<0$,
then the Heegner cycle $\gamma_{D_0,D_1,\rho}$ is
anti-invariant under the action of complex conjugation on $X_0^+(N)(\CC)$.
In particular the homology class $[\gamma_{D_0,D_1,\rho}]$ represented by
the cycle $\gamma_{D_0,D_1,\rho}$
in fact lies in the $-1$-eigenspace $H_1(X_0^+(N)(\CC),\ZZ)^{-}$.
Following Gross-Kohnen-Zagier~\cite[p. 559]{gross:gkz} we may define
an element $e\in H_1(E(\CC),\ZZ)^{-}$ such that
\begin{equation}
[\gamma(D_0,D_1,r_0r_1)]_E=I_{D_0,r_0,E}I_{D_1,r_1,E}\,e_E,
\end{equation}
where $[\gamma(D_0,D_1,r_0r_1)]_E$ is the canonical image in $H_1(E(\CC),\ZZ)^{-}$
of the homology class $[\gamma(D_0,D_1,r_0r_1)]$,
and as above $I_{D_i,r_i,E}$ denotes the index of the trace $y_{D_i,r_i,E}$ in $E(\QQ)$,
for all pairs $(D_i,r_i)$ with $D<0$ that satisfy the Heegner condition.
It is well-known that the index $n_E$
of the subgroup generated by $e_E$ in $H_1(E(\CC),\ZZ)^{-}$
is uniquely defined by the above condition.

Ogg~\cite{ogg:real} describes the real locus $(S/w_m)(\RR)$ of
quotients $S/w_m$ of Shimura curves $S$,
attached to Eichler orders $\mathcal{O}$ of indefinite quaternion algebras over $\QQ$,
in terms of embeddings of $\QQ(\sqrt{m})$ into $\mathcal{O}$.
In particular,
from his work it is known that
the number $\nu_N$ of connected components of $X_0^+(N)(\RR)$
is given by the formula
\begin{displaymath}
\nu_N=
\left\{
\begin{array}{ll}
\frac{h(4N)+h(N)}{2},&\textrm{if $N\equiv 1\pmod{4}$}\\
\\
\frac{h(4N)+1}{2},&\textrm{otherwise.}\\
\end{array}
\right.
\end{displaymath}
Moreover,
as shown in~\cite{ccb:thesis}
it is possible to describe explicitly
the connected components of $X_0^+(N)(\RR)$
as a sum of ``weighted'' Heegner cycles
over discriminants $\Delta > 0$ such that $N|\Delta$ and $\Delta|4N$,
in analogy with the fixed points of the Fricke involution
(cf. Gross~\cite{gross:primelevel}).

%
\subsection{The conjecture}\label{sub:dis}
As above,
let $E$ be an elliptic curve of rank one over $\QQ$,
and let $I_\varphi$  be the index of
the group generated by the Heegner points,
i.e. the greatest common divisor of the indexes $I_{D,r,\varphi}$
for all pairs $(D,r)$ that satisfy the Heegner condition
with fundamental $D<0$.
From now on assume that $N_E$ is prime.
In particular $E$ is alone in its isogeny class,
so we may write  $I_E$  instead of $I_\varphi$.


\begin{conjecture}\label{con:main}
If $I_E>1$ then either the number $\nu_{N_E}$ of connected components of
the real locus $X_0^+(N_E)(\RR)$ is $\nu_{N_E}>1$
or the Tate-\v Safarevi\v c group $\Sha(E)$ of $E$ is non-trivial.
\end{conjecture}

There are some curves $E$ in the range of our computations
that have $\nu_{N_E}>1$ but have index $I_E=1$.
So knowing $\nu_N$ is not enough in order to predict when $I_E > 1$.
We sketch,
in a rather impressionistic style,
some ideas that hopefully will lead to
a more aesthetically pleasing version of the conjecture
as follows.
As shown by Gross-Harris~\cite[pp. 164--165]{gross:real},
given any complete, non-singular, geometrically connected curve
defined over $\RR$
the number $\nu$ of connected components of $X(\RR)$
may be recovered from the homology group $H_1(X(\CC),\FF_2)$,
regarded as a symplectic $\FF_2$-vector space
with involution $\tau$
induced by complex conjugation acting on $X(\CC)$.
In fact they prove that
\begin{equation*}
\nu=g+1-\textrm{rank}(H)
\end{equation*}
where $g$ is the genus of $X$,
and $H$ is the $g$ by $g$ symmetric matrix defined by
\begin{equation*}
[\tau]_\beta=\mtwo{I_g}{H}{0}{I_g},
\end{equation*}
where $\beta$ is a suitable symplectic basis for $H_1(X(\CC),\FF_2)$.
So our conjecture may be expressed in homological terms.
It is hoped that a more refined version of our conjecture may be
meaningfully stated in terms of a finer homological invariant,
perhaps associated to
the modular parametrisation $X_0^+(N)\longrightarrow E$ over $\QQ_p$
for each prime $p$,
with special attention to the primes $p=N,\infty$;
maybe there is some kind of ``product formula'' for $n_E$
in which $\nu_N$ is just a very crude approximation to
the contribution from $p=\infty$.
Such formula might lead to a more natural form of Equation~\ref{eqn:gkz},
if we consider that
the Tate-\v Safarevi\v c group $\Sha(E)$ is a subgroup of
the cohomology group $H^1(G_\QQ, E)$
determined by local conditions.

Table~\ref{tbl:indexes01} and Table~\ref{tbl:indexes02}
(below) were computed as follows.
For each elliptic curve $E$ of rank one over $\QQ$
and prime conductor $N_E < 129999$,
we computed the greatest common divisor $d$ of
the indexes $I_{D,E}$,
for each pair $(D, r)$ that satisfies the Heegner condition
with $D<0$ fundamental and $|D|\leq 163$.
Such $d$ is likely to be the index $I_E$
of the group generated by
all the traces $y_{D,E}$ in $E(\QQ)$ in the range $N_E < 129999$.
All our elliptic curve data comes from Cremona's Tables~\cite{cremona:onlinetables},
and we stick to the notation used there.
\begin{table}
\ssp
\begin{center}
\caption{Nontrivial indexes $I_E$ for prime $N_E\leq 84701$.
\label{tbl:indexes01}}
\end{center}%
\vspace{-.3in}%
$$
\begin{array}{llll}
E &	I_{E}&	\nu_{N_E}& \Sha(E) \\
\vspace{-2ex}\\
{\bf 359 A}&		2&	2&	1\\
{\bf 359 B}&		2&	2&	1\\
{\bf 997 A}&		2&	2&	1\\
{\bf 3797 A}&		2&	2&	1\\
{\bf 4159 A}&		2&	2&	1\\ 
\vspace{-2ex}\\
{\bf 4159 B}&		2&	2&	1\\
{\bf 6373 A}&		2&	2&	1\\
{\bf 8069 A}&		2&	3&	1\\
{\bf 8597 A}&		2&	6&	1\\
{\bf 9829 A}&		2&	10&	1\\ 
\vspace{-2ex}\\
{\bf 13723 A}&		2&	2&	1\\
{\bf 17299 A}&		2&	2&	1\\
{\bf 17573 A}&		2&	2&	1\\
{\bf 18097 A}&		2&	3&	1\\
{\bf 18397 A}&		2&	2&	1\\ 
\vspace{-2ex}\\
{\bf 20323 A}&		2&	2&	1\\
{\bf 21283 A}&		2&	2&	1\\
{\bf 23957 A}&		2&	6&	1\\
{\bf 24251 A}&		2&	5&	1\\
{\bf 26083 A}&		2&	2&	1\\ 
\vspace{-2ex}\\
{\bf 28621 A}&		2&	2&	1\\
{\bf 28927 A}&		2&	2&	1\\
{\bf 29101 A}&		2&	2&	1\\
{\bf 29501 A}&		2&	2&	1\\
{\bf 31039 A}&		2&	2&	1\\ 
\vspace{-2ex}\\
{\bf 31319 A}&		2&	2&	1\\
{\bf 33629 A}&		2&	2&	1\\
{\bf 34613 A}&		2&	2&	1\\
{\bf 34721 A}&		2&	3&	1\\
{\bf 35083 B}&		4&	1&	4\\ 
\vspace{-2ex}\\
{\bf 35401 A}&		2&	3&	1\\
{\bf 35533 A}&		2&	2&	1\\
{\bf 36479 A}&		2&	11&	1\\
{\bf 36781 A}&		2&	2&	1\\
{\bf 36781 B}&		2&	2&	1\\ 
\end{array}\qquad
\begin{array}{llll}
E &	I_{E}&	\nu_{N_E}& \Sha(E) \\
\vspace{-2ex}\\
{\bf 39133 A}&		2&	2&	1\\
{\bf 39133 B}&		2&	2&	1\\
{\bf 39301 A}&		2&	14&	1\\
{\bf 40237 A}&		2&	2&	1\\
{\bf 45979 A}&		4&	2&	4\\ 
\vspace{-2ex}\\
{\bf 47143 A}&		2&	2&	1\\
{\bf 47309 A}&		2&	2&	1\\
{\bf 48731 A}&		4&	1&	4\\
{\bf 50329 A}&		2&	3&	1\\
{\bf 51437 A}&		2&	6&	1\\ 
\vspace{-2ex}\\
{\bf 52237 A}&		2&	2&	1\\
{\bf 55837 A}&		2&	14&	1\\
{\bf 59243 A}&		2&	2&	1\\
{\bf 61909 A}&		2&	6&	1\\
{\bf 62191 A}&		2&	5&	1\\ 
\vspace{-2ex}\\
{\bf 63149 A}&		2&	2&	1\\
{\bf 65789 A}&		2&	2&	1\\
{\bf 66109 A}&		2&	2&	1\\
{\bf 66109 B}&		2&	2&	1\\
{\bf 67427 A}&		2&	5&	1\\ 
\vspace{-2ex}\\
{\bf 68489 B}&		2&	3&	1\\
{\bf 69677 A}&		2&	2&	1\\
{\bf 72053 A}&		2&	2&	1\\
{\bf 73709 A}&		2&	2&	1\\
{\bf 74411 A}&		2&	2&	1\\ 
\vspace{-2ex}\\
{\bf 74713 A}&		4&	3&	4\\
{\bf 74797 A}&		2&	2&	1\\
{\bf 77849 A}&		2&	3&	1\\
{\bf 78277 A}&		2&	2&	1\\
{\bf 78919 A}&		2&	2&	1\\ 
\vspace{-2ex}\\
{\bf 81163 B}&		2&	2&	1\\
{\bf 81349 A}&		2&	2&	1\\
{\bf 82301 A}&		2&	2&	1\\
{\bf 84653 A}&		2&	2&	1\\
{\bf 84701 A}&		2&	3&	1\\ 
\end{array}
$$
\end{table}

\begin{table}
\ssp
\begin{center}
\caption{Nontrivial indexes $I_E$ for prime $85837\leq N<129999$.
\label{tbl:indexes02}}
\end{center}%
\vspace{-.3in}%
$$
\begin{array}{llll}
E &	I_{E}&	\nu_{N_E}& \Sha(E) \\
\vspace{-2ex}\\
{\bf 85837 A}&		2&	2&	1\\
{\bf 87013 A}&		2&	3&	1\\
{\bf 90001 B}&		2&	87&	1\\
{\bf 90001 C}&		2&	87&	1\\
{\bf 90001 D}&		2&	87&	1\\ 
\vspace{-2ex}\\
{\bf 91381 A}&		2&	2&	1\\
{\bf 92419 A}&		4&	1&	4\\
{\bf 101771 A}&		2&	2&	1\\
{\bf 101879 A}&		2&	2&	1\\
{\bf 102061 B}&		2&	6&	1\\ 
\vspace{-2ex}\\
{\bf 103811 A}&		2&	2&	1\\
{\bf 104239 A}&		4&	14&	4\\
{\bf 104239 B}&		4&	14&	4\\
{\bf 105143 A}&		2&	2&	1\\
{\bf 105401 A}&		2&	3&	1\\ 
\vspace{-2ex}\\
{\bf 105541 A}&		2&	2&	1\\
{\bf 106277 A}&		2&	14&	1\\
{\bf 106949 A}&		2&	2&	1\\
{\bf 106979 A}&		4&	1&	4\\
{\bf 107981 A}&		2&	2&	1\\ 
\end{array}\qquad
\begin{array}{llll}
E &	I_{E}&	\nu_{N_E}& \Sha(E) \\
\vspace{-2ex}\\
{\bf 108971 A}&		2&	2&	1\\
{\bf 113933 A}&		2&	2&	1\\
{\bf 118673 A}&		2&	3&	1\\
{\bf 119689 A}&		2&	3&	1\\
{\bf 119701 A}&		2&	3&	1\\ 
\vspace{-2ex}\\
{\bf 119773 A}&		2&	2&	1\\
{\bf 123791 A}&		2&	2&	1\\
{\bf 124213 A}&		2&	2&	1\\
{\bf 126683 A}&		2&	2&	1\\
{\bf 127669 A}&		2&	2&	1\\ 
\vspace{-2ex}\\
{\bf 129277 A}&		2&	2&	1\\
{\bf 129853 A}&		2&	2&	1\\
\\
\\
\\
\vspace{-2ex}\\
\\
\\
\\
\\
\\
\end{array}
$$
\end{table}

%
%
%
\bibliographystyle{amsplain}
\bibliography{biblio}
\end{document}